\documentclass{amsart}
\usepackage{amssymb, latexsym}

\DeclareMathOperator{\Alt}{Alt}

\DeclareMathOperator{\M}{\mathcal{M}}
\DeclareMathOperator{\N}{\mathcal{N}}

\DeclareMathOperator{\Aut}{Aut}

\theoremstyle{plain}
\newtheorem{theorem}{Theorem}
\newtheorem{corollary}{Corollary}
\newtheorem{lemma}{Lemma}

\theoremstyle{definition}

\begin{document}
\title[Subspaces of $7\times 7$ matrices]
{Subspaces of $7\times 7$ skew-symmetric\\ 
matrices related to the group $G_2$} 

\author{Rod Gow}
\address{School of Mathematical Sciences\\
University College\\
Belfield, Dublin 4\\
Ireland}
\email{rod.gow@ucd.ie}
\keywords{rank, constant rank subspace, octonions, skew-symmetric matrix, automorphisms}

\subjclass{15A03, 15A33}
\begin{abstract} 
Let $K$ be a field of characteristic different from 2 and let $C$ be an octonion algebra over $K$. We show that there is a seven-dimensional subspace
of $7\times 7$ skew-symmetric matrices over $K$ which is invariant under the automorphism group of $C$. This subspace consists of elements of rank 6 when $C$ is a division algebra, and elements of rank 4 and 6 when $C$ is a split algebra.
In the latter case, the automorphism group is the exceptional group $G_2(K)$.

\end{abstract}
\maketitle
\section{Introduction}
\noindent  Let $K$ be a field of characteristic different from 2 and let $K^*$ denote the subset of non-zero elements
in $K$. Let
$n$  be a positive integer and let $M_n(K)$ denote the algebra of $n\times n$ matrices with entries in $K$. Let $r$ be an integer satisfying $1\le r\le n$. We say that a non-zero subspace $\mathcal{M}$ of $M_n(K)$
is a \emph{constant rank} $r$ subspace if each non-zero element of $\mathcal{M}$ has rank $r$. 

Such subspaces have been the subject of much research, in linear algebra, algebraic geometry and differential topology, for example. The references of
\cite{LY} provide information on earlier work about constant rank subspaces 

Lam and Yiu provided an example of a seven-dimensional constant rank 6 subspace of
skew-symmetric matrices in $M_7(\mathbb{R})$, \cite{LY}, Proposition 10. This subspace is obtained from 
a seven-dimensional constant rank 8 subspace of skew-symmetric matrices in $M_8(\mathbb{R})$, derived from the regular representation of the division algebra of real octonions on itself.

 The purpose of this paper to generalize the construction of Lam and Yiu, putting it into the context of octonion
algebras and their automorphisms. When there is an octonion division algebra over $K$, there is a corresponding 
seven-dimensional constant rank 6 subspace of skew-symmetric matrices in $M_7(K)$ which is invariant under the action of the automorphism group
of the division algebra. (Note that there may be non-isomorphic octonion division algebras for appropriate $K$.) 

For the split octonion algebra, which exists for all fields $K$, there is a 
seven-dimensional subspace of skew-symmetric matrices in $M_7(K)$ in which all non-zero elements have rank 6 or 4. This subspace
is also invariant under the action of the automorphism group of the algebra, which in this case is the split group  $G_2(K)$.

Furthermore, when there is an octonion division algebra over $K$, we show that there is  
a 14-dimensional subspace of skew-symmetric matrices in $M_7(K)$ containing no elements of rank 2. This subspace is also
invariant under the action of the automorphism group of division algebra. We note here that
14 is the maximum dimension for such a subspace of skew-symmetric matrices. 

We phrase our proofs in terms of alternating bilinear forms, but the correspondence between these and skew-symmetric
matrices should enable the reader to reinterpret our results in terms of matrices. We also make extensive use of the notation and theory of Chapters 1 and 2 of the book \cite{SV}, although much of our work is reasonably self-contained.


\section{Properties of octonion algebras and their automorphisms}

\noindent We shall work in the context of an octonion algebra $C$ over $K$.
\smallskip
$C$ is an eight-dimensional non-associative algebra over $K$
with multiplicative identity $e$. It is equipped with a non-degenerate quadratic form $N$ which permits composition, meaning that 
\[
 N(xy)=N(x)N(y)
\]
for all $x$ and $y$ in $C$. We call $N$ the \emph{norm} on $C$. It satisfies $N(e)=1$. We let $\langle\,,\,\rangle$ denote the associated polarizing symmetric bilinear form. 

Let $x$ be an element of $C$. We define the conjugate $\overline{x}$ of $x$ by
\[
  \overline{x}=\langle x,e\rangle e-x
\]
The conjugation operator is an anti-involution on 
 $C$ satisfying
\[
 x\overline{x}=\overline{x}x=N(x)e.
\]
Provided that $N(x)\ne 0$, $x$ has an inverse $x^{-1}$ given by
\[
 x^{-1}=\frac{1}{N(x)}\overline{x}.
\]

We say that the octonion $x$ is \emph{pure} if 
\[
\overline{x}=-x.
\] 
This is equivalent to saying that $\langle x,e\rangle=0$. The pure octonions form a seven-dimensional subspace of
$C$, which we shall denote by $C_0$. The restriction of $N$ to $C_0$ is clearly non-degenerate.

Let $x$ be an octonion. We define a $K$-linear transformation $L_x$ from $C$ into itself by
\[
 L_x(y)=xy
\]
for all octonions $y$. Properties of octonions imply that $L_x$ is invertible if and only $x$ is invertible. The following facts concerning
the kernel and image of $L_x$ when $x$ is not invertible are not new, but we supply a proof
for the sake of completeness.

\begin{lemma} \label{kernelandimage} Let $x$ be a non-zero non-invertible octonion and let $I_x$ denote the image of 
$L_x$. Then $\dim I_x=4$. Furthermore, 
\[
 \dim\,(I_x\cap\,C_0)=\dim\,(\ker L_x\cap\, C_0)=3.
\]
\end{lemma}

\begin{proof} An element of $I_x$ has the form $xy$. Since $N(xy)=N(x)N(y)$ and $N(x)=0$, $xy$ is isotropic
with respect to $N$. It follows that $I_x$ is totally isotropic and since $N$ is non-degenerate on $C$, we have
\[
 \dim I_x\le 4.
\]
Similarly, let $z$ be an element of $\ker L_x$. Since $xz=0$, $z$ cannot be a unit and is hence isotropic. It follows that
$\ker L_x$ is also totally isotropic, and we similarly deduce that $\dim \ker L_x\le 4$. As we have 
\[
 \dim I_x+\dim \ker L_x=8
\]
by the rank-nullity theorem, we see that
\[
  \dim I_x=\dim \ker L_x=4.
\]
Finally, since $I_x\cap\, C_0$ and $\ker L_x\cap C_0$
are totally isotropic subspaces of the seven-dimensional non-degenerate space $C_0$, they both have
dimension at most 3. On the other hand, the intersections have dimension at least 3, since $C_0$ has codimension 1 in $C$.
\end{proof}

An \emph{automorphism} $\sigma$ of $C$ is an injective $K$-linear transformation from $C$ onto itself satisfying
\[
 \sigma(xy)=\sigma(x)\sigma(y)
\]
for all $x$ and $y$. The automorphisms of $C$ form a group, which we shall denote by $\Aut(C)$. When $C$ is a split algebra (in other words, when $N$ is isotropic),
$\Aut(C)$ is the split group $G_2(K)$ over $K$. When $C$ is a division algebra,
the isomorphism type of the automorphism group depends on the isometry class of $N$, but we can say somewhat imprecisely
that the group is of $G_2$ type.

Each automorphism $\sigma$ of $C$ fixes $e$ and is an isometry of $N$; see, for example, Corollary 1.2.4 of \cite{SV}. It follows that $\sigma$ maps the subspace of pure octonions onto itself, since it is the orthogonal complement of $e$.

\section{Construction of a subspace of alternating bilinear forms}

\noindent Let $x$ be a pure octonion. We define a bilinear form
\[
 F_x: C\times C\to K
\]
by setting 
\[
 F_x(y,z)=\langle xy,z\rangle.
\]

We now prove some simple properties of the bilinear form $F_x$.

\begin{lemma} \label{alternating} The bilinear form $F_x$ is alternating and its radical is $\ker L_x$. Consequently,
$F_x$ has rank $8$ or $4$. Furthermore,
\[
 F_x(y,z)=-\langle x, y\overline{z} \rangle.
\]

\end{lemma}

\begin{proof}
 By Lemma 1.3.2 of \cite{SV}, we have
\[
F_x(y,z)= \langle xy,z\rangle=\langle y,\overline{x}z\rangle.
\]
Now $\overline{x}=-x$ since $x$ is pure, and thus
\[
 F_x(y,z)=-\langle y,xz\rangle=-\langle xz,y \rangle=-F_x(z,y).
\]
This implies that $F_x$ is alternating. 

\smallskip

Suppose now that $y$ is in the radical of $F_x$. Then we have 
\[
 \langle xy,z\rangle=0
\]
for all $z\in C$. Since the polarizing form is non-degenerate, $xy=0$ and hence $y\in \ker L_x$. Conversely,
any element of $\ker L_x$ is in the radical of $F_x$, This implies that the radical is $\ker L_x$. 

The statement about the rank of $F_x$ now follows from Lemma \ref{kernelandimage}. The final formula of the lemma follows from formula 1.12 of Lemma 1.3.2 of \cite{SV}. 

\end{proof}

Let $\Alt(C)$ denote the vector space of all alternating bilinear forms $C\times C\to K$. We define a function
\[
 \epsilon: C_0\to \Alt(C)
\]
by setting 
\[
 \epsilon(x)=F_x.
\]
It is elementary to see that $\epsilon$ is linear and injective. Let $\mathcal{N}$ denote the image of $\epsilon$. Lemma \ref{alternating} implies that $\mathcal{N}$ is a seven-dimensional subspace of $\Alt(C)$ whose non-zero elements have rank 8 or 4. Clearly, $\mathcal{N}$ is a constant rank 8 subspace precisely when $C$ is a division algebra.

Given a pure octonion $x$ and corresponding bilinear form $F_x$, let $f_x$ denote the restriction of $F_x$ to
$C_0\times C_0$. It is clear that $f_x$ is non-zero if $x$ is non-zero. Thus if $\M$ denotes the subset of all such $f_x$,
$\M$ is a seven-dimensional subspace of $\Alt(C_0)$. 

To find the rank of $f_x$, we make use of the following general principle, which must be well known. We include a proof for the  convenience of the reader.

\begin{lemma} \label{restriction} Let $V$ be a finite dimensional vector space over $K$ and let $U$ be a subspace of codimension $1$ in $V$. Let $F$ be an alternating bilinear form defined on $V\times V$ and let $f$ be its restriction to $U\times U$. Let $2r$ be the rank of $F$. Then $f$ has rank $2r$ or $2r-2$, and the latter rank occurs precisely when the radical of $F$ is contained in $U$.
\end{lemma}

\begin{proof} Suppose first that $F$ is non-degenerate. Then $\dim V=2r$ and $\dim U=2r-1$. Since $U$ has odd dimension
 and $f$ must have even rank, $f$ has rank at most $2r-2$. Now the radical of $f$ is $U\cap U^\perp$. Since $F$ is non-degenerate, we have
\[
 \dim V=2r=\dim U+\dim U^\perp=2r-1+\dim U^\perp.
\]
Thus $\dim U^\perp=1$ and we see that the radical of $f$ is one-dimensional. This implies that $f$ has rank
$2r-2$. 

We now consider the general case. Let $R$ be the radical of $F$. Suppose that $R$ is contained in $U$. Let $W$ be a complement of $R$ in $V$ and let $G$ be the restriction of $F$ to $W\times W$.  Then $G$ is a non-degenerate
alternating bilinear form of rank $2r$ which determines $F$. We also have $U=R\oplus U\cap W$ and $U\cap W$ has codimension 1 in $W$. The restriction of $G$ to $U\cap W\times U\cap W$ has rank $2r-2$ by the argument above. This implies that
$f$ has rank $2r-2$ also. 

Finally, suppose that $R$ is not contained in $U$. Let $S$ be a one-dimensional subspace of $R$ not contained in $U$. Then we have
\[
 V=S\perp U
\]
and, since $F$ is determined by its restriction to $U\times U$, $F$ and $F$ have the same rank, $2r$, in this case.
\end{proof}

We apply this lemma to the forms $f_x$ constructed above.

\begin{lemma}  Let $x$ be a non-zero pure octonion. Then $f_x$ has rank $6$ or $4$. Moreover $f_x$ has rank $4$ precisely when $x$ is not invertible. Thus $\M$ is a constant rank $6$ subspace precisely when $C$ is a division algebra.
\end{lemma}

\begin{proof} Suppose that $x$ is invertible. Then
 $F_x$ has rank 8, and it follows from Lemma \ref{restriction} that $f_x$ has rank 6. Suppose on the other hand, that
$x$ is not invertible. It follows from Lemma \ref{kernelandimage} that $F_x$ has rank 4 and its radical is $\ker L_x$. We also know from Lemma \ref{kernelandimage} that $\ker L_x$ is not contained in $C_0$. It follows from Lemma \ref{restriction} that $f_x$ also has rank $4$.
\end{proof}

Let $\sigma$ be an element of $\Aut(C)$ and let $f$ be an element of $\Alt(C)$ or $\Alt(C_0)$. We define the alternating bilinear form $\sigma(f)$ by
\[
 \sigma(f)(x,y)=f(\sigma^{-1}(x),\sigma^{-1}(y)).
\]
The formation of $\sigma(f)$ from $f$ defines a linear action (representation) of $\Aut(C)$ on the spaces $\Alt(C)$ and
$\Alt(C_0)$. 

We now show that the subspaces $\M$ and $\N$ are invariant under $\Aut(C)$. 

\begin{lemma} Let $x$ be a pure octonion and let $F_x$ be the corresponding element of $\M$. Then 
 \[
  \sigma(f_x)=f_{\sigma(x)}.
 \]
Hence $\M$ is $\Aut(C)$-invariant. An identical conclusion holds for $F_x$ and $\N$. 

\end{lemma}

\begin{proof} By definition,
 \[
  \sigma(f_x)(y,z)=f_x(\sigma^{-1}(y),\sigma^{-1}(z))
 \]
for all $y$ and $z$ in $C_0$. Thus
\[
  \sigma(f_x)(y,z)=\langle x\sigma^{-1}(y),\sigma^{-1}(z)\rangle.
 \]
However, $\sigma$ is an isometry of the polarizing form and hence
\[
 \langle x\sigma^{-1}(y),\sigma^{-1}(z)\rangle=\langle \sigma(x)y,z\rangle.
\]
This implies that
\[
  \sigma(f_x)=f_{\sigma(x)},
 \]
as required. The rest of the proof is the same.

\end{proof}

We can now summarize our findings.

\begin{theorem} \label{subspace} Let $C$ be an octonion algebra over the field $K$ and let $\Aut(C)$ denote the automorphism group of $C$. Let $C_0$ denote the seven-dimensional subspace of pure octonions. Then there is a seven-dimensional subspace $\M$ of alternating bilinear forms defined on $C_0\times C_0$ which is invariant under $\Aut(C)$. All the non-zero elements of $\M$ have rank $6$ when $C$ is a division algebra. 

When $C$ is the split octonion algebra, $\Aut(C)$
 is then the split group $G_2(K)$ and the non-zero elements of $\M$ have rank $6$ or $4$ (and both ranks occur).
\end{theorem}

The action of $\Aut(C)$ on $\M$ can be better appreciated when we take into account the group's action on pure octonions.
The following result must be well known, but for want of a convenient explicit reference, we provide a proof.

\begin{lemma} \label{transitivity}
 Let $x$ and $y$ be non-zero pure octonions with $N(x)=N(y)$. Then there exists $\sigma\in\Aut(C)$ with $\sigma(x)=y$.
\end{lemma}

\begin{proof} Suppose first that $N(x)\ne 0$. Then the result follows from Corollary 1.7.5 of \cite{SV}.
 Suppose next that $N(x)=0$. By elementary properties of symmetric bilinear forms, working in the subspace
$C_0$ we may find a pure octonion $z$ with 
\[
 \langle x,z\rangle =1, \quad N(z)=0.
\]
We set
\[
 a=x+z,\quad b=x-z.
\]
Since $N(x)=N(z)=0$, we have 
\[
 \langle x,x\rangle= \langle z,z\rangle=0
\]
and then we find that
\[
  \langle a,b\rangle=0.
\]
Since $K$ has characteristic different from 2, we also obtain
\[
 N(a)=1,\quad N(b)=-1.
\]
We note also that $a$ and $b$ are both pure.

We may likewise construct a pure octonion $w$, say, with
\[
 \langle y,w\rangle =1, \quad N(w)=0.
\]
We then set
\[
 a'=y+w,\quad b'=y-w
\]
and calculate as above that
\[
 N(a')=1,\quad N(b')=-1,\quad \langle a',b'\rangle=0
\]
It follows from Corollary 1.7.5 of \cite{SV} that
there is an element $\sigma$ of $\Aut(C)$ satisfying
\[
 \sigma(a)=a',\quad \sigma(b)=b'.
\]
But now as
\[
 x=\frac{a+b}{2},\quad y=\frac{a'+b'}{2}
\]
and $\sigma$ is $K$-linear, we have
\[
 \sigma(x)=y,
\]
as required.

\end{proof}

\begin{corollary} \label{isotropicpoints}
 Suppose that $C$ is a split octonion algebra over $K$. Then $\Aut(C)$ acts transitively on the non-zero non-invertible pure
octonions.
\end{corollary}

Rather than looking at the $\Aut(C)$ action on individual invertible pure octonions, we look instead at the action on the one-dimensional subspaces they generate. Here, the orbits depend on the square classes in $K^*$, and consequently, we restrict attention to some important fields. 

\begin{lemma}
$\Aut(C)$ has exactly one orbit on the one-dimensional subspaces of $C_0$ generated by invertible octonions when
$K$ is algebraically closed, or when $K=\mathbb{R}$ and $C$ is a division algebra.
$\Aut(C)$ has exactly two orbits on the one-dimensional subspaces of $C_0$ generated by invertible octonions when
$K$ is finite, or when $K=\mathbb{R}$ and $C$ is a split algebra. 
\end{lemma}

\begin{proof} When $K$ is algebraically closed, each non-zero scalar is a square. Consequently, for a given non-zero octonion $x$, all elements of $K^*$ are expressible as $N(\lambda x)$ as $\lambda$ runs over $K^*$. We obtain
 the desired transitivity conclusion from Lemma \ref{transitivity}. Suppose next that $K=\mathbb{R}$ and $C$ is a division algebra. Then $N$ is positive definite, so that $N(x)$ is a positive real number whenever $x$ is non-zero. The result follows in this case as each positive real number is a square. Finally, when $K$ is finite, $C$ is automatically split
and since there are two square classes in $K^*$ in this case, the result follows by the previous reasoning.
\end{proof}

\begin{corollary} \label{orbits} $\Aut(C)$ has exactly one orbit on the bilinear forms of rank $4$ in $\M$ (such forms
exist only when $C$ is a split algebra).
$\Aut(C)$ has exactly one orbit on the one-dimensional subspaces of $\M$ generated by bilinear forms of rank $6$ when
$K$ is algebraically closed. When $K=\mathbb{R}$ and $C$ is a division algebra, $\Aut(C)$ also has exactly one orbit on all the one-dimensional subspaces of $\M$.
$\Aut(C)$ has exactly two orbits on the one-dimensional subspaces of $\M$ generated by bilinear forms of rank $6$ when
$K$ is finite, or when $K=\mathbb{R}$ and $C$ is a split algebra. 
\end{corollary}

We remark that when $K=\mathbb{R}$, $\Aut(C)$ is the compact real Lie group of type $G_2$ if $C$ is a division algebra, and the non-compact real Lie group of type $G_2$ if $C$ is split. These groups are denoted by $G_2^c$ and $G_2^*$ in the instructive article \cite{A}. Some of the geometry associated to these two groups can be explained in terms of their actions considered here.

\section{A $14$-dimensional subspace of alternating bilinear forms}
\noindent
Suppose that $C$ is a division algebra over $K$. Then the subspace $\M$ of $\Alt(C_0)$ described in Theorem \ref{subspace}
has the unusual property that it does not vanish identically on any two-dimensional subspace of $C_0$, as we show below.

\begin{lemma} \label{scalarmultiple} Suppose that $C$ is a division algebra over $K$. Let $y$ and $z$ be non-zero pure octonions. Suppose also that
 \[
  f_x(y,z)=0
 \]
for all elements $x$ of $C_0$. Then $y=\lambda z$ for some element $\lambda$ of $K$.
\end{lemma}

\begin{proof}
 Suppose that
 \[
  f_x(y,z)=0
 \]
for all elements $x$ of $C_0$. It follows from Lemma \ref{alternating} that
\[
 \langle x, y\overline{z}\rangle=0
\]
and thus $y\overline{z}$ is orthogonal to $C_0$ with respect to the polarizing form. But the orthogonal complement of
$C_0$ is the subspace spanned by the identity element $e$. Hence $y\overline{z}=\mu e$ for some scalar $\mu$. Multiplying on the right by $z$, we see that $N(z)y=\mu z$ and the result follows.
\end{proof}

\begin{corollary} \label{nonvanishing}
 There is no two-dimensional subspace of $C_0$ on which the elements of $\M$ all vanish.
\end{corollary}

Let $C_0\wedge C_0$ denote the exterior square of $C_0$. We say that an element $z$ of $C_0\wedge C_0$ is
\emph{decomposable} if we have $z=x\wedge y$ for suitable
$x$ and $y$ in $C_0$.

Let $f$ be an element of $\Alt(C_0)$ and let $\{e_1, \ldots, e_7\}$
be a basis of $C_0$. We may define a linear form $f^*$ on
$C_0\wedge C_0$ by setting
\[
f^*(e_i\wedge e_j)=f(e_i,e_j)
\]
and extending to all of $C_0\wedge C_0$ by linearity.
We note then that
\[
f^*(x\wedge y)=f(x,y)
\]
for any $x$ and $y$.

For the sake of simplicity, we write $f_i$ in place of $f_{e_i}$. We now define
a linear transformation 
\[
\omega: C_0\wedge C_0\to K^7
\]
by
\[
\omega(z)=(f_1^*(z), \ldots, f_7^*(z))
\]
for all $z\in C_0\wedge C_0$.

\begin{theorem} \label{kernel} With the notation previously introduced,
$\omega$ is surjective and the kernel of $\omega$ is a subspace
of codimension $7$ in $C_0\wedge C_0$ which contains no non--zero
decomposable elements and is invariant under $\Aut(C)$.
\end{theorem}

\begin{proof} We first show that $\omega$ is surjective.
Let $v_i$ be the standard basis vector of $K^7$ whose single
non--zero component is 1 occurring in the $i$--th position.
We show that $v_i$ is in the image of $\omega$. Now as $C_0$ is seven-dimensional,
given any six elements of $\Alt(C_0)$, there is a two-dimensional subspace of $C_0$ on which they all vanish.
Thus there is a two--dimensional subspace $U_i$, say, of
$C_0$ that is isotropic for the six forms $f_1$, \dots,
$f_{i-1}$, $f_{i+1}$, \dots, $f_7$ but is not isotropic
for $f_i$ (this last statement follows from Corollary \ref{nonvanishing}). Let $x_i$, $y_i$ be basis vectors for $U_i$
with $f_i(x_i,y_i)=1$. Then we have
\[
\omega(x_i\wedge y_i)=v_i,
\]
as required. 

We next show that $\ker \omega$  contains no non--zero
decomposable element. For suppose that
\[
\omega(x\wedge y)=0.
\]
Then we have
\[
f_1(x,y)=\cdots =f_7(x,y)=0.
\]
Since the $f_i$ do not simultaneously vanish on any
two-dimensional subspace of $C_0$, by Corollary \ref{nonvanishing}, $x$ and $y$ must be linearly
dependent and hence $x\wedge y=0$.

Finally, suppose that $z\in \ker \omega$. Then 
\[
f_1^*(z)=\cdots =f_7^*(z)=0.
\]
Since $f_1$, \dots, $f_7$ are a basis of $\M$, we have $f^*(z)=0$ for all $f\in\M$. Now as
$\M$ is $\Aut(C)$-invariant and $\Aut(C)$ acts on $C_0\wedge C_0$ by means of $x\wedge y\to \sigma(x)\wedge \sigma(y)$,
we have
\[
 \sigma(f)^*(z)=0=f^*(\sigma^{-1}(z))
\]
for all $f\in\M$ and all $\sigma\in\Aut(C)$. This implies that $\ker\, \omega$ is invariant under $\Aut(C)$.

\end{proof}

We conclude this section by explaining the connection between Theorem \ref{kernel} and the earlier work. Let $V$ be a finite dimensional vector space over $K$. It is 
well known that there is a linear isomorphism between $V\wedge V$ 
and $\Alt(V^*)$, where $V^*$ is the dual space of $V$.
The isomorphism is defined in the following way. 

Let
$\{v_1,\dots ,v_n\}$ be a basis of $V$ and let 
$$
z=\sum_{1\le i<j\le n}a_{ij}v_i\wedge v_j
$$
be any element of $V\wedge V$. Define 
$$
\varepsilon_z : V^*\times V^*\longrightarrow K
$$
by 
$$
\varepsilon_z(\theta,\phi ) = \sum
a_{ij}\left(\theta(v_i)\phi(v_j)-\theta(v_j)\phi(v_i)\right).
$$
for all $\theta$ and $\phi$ in $V^*$.

It is straightforward to verify that $\varepsilon_z$ is in 
$\Alt(V^*)$ and the mapping $z\to\varepsilon_z$ is
an isomorphism between $V\wedge V$ and $\Alt(V^*)$.

The rank of $\varepsilon_z$ is the
dimension of the subspace of $V$ associated with $z$. 
(Recall that the subspace  associated
to $z$ is the smallest subspace
$U$ of $V$ such that $z\in U\wedge U$.) Thus
$\varepsilon_z$ has rank $2$ precisely when $z$ is non--zero
and decomposable.

Now $C_0$ and $C_0^*$ are isomorphic as $\Aut(C)$-modules, the isomorphism being provided by the bilinear form
$ \langle ,\rangle$. Consequently, $C_0\wedge C_0$ and $\Alt(C_0)$ are also isomorphic as $\Aut(C)$-modules.

We can thus reinterpret Theorem \ref{kernel} in the following terms.

\begin{theorem} \label{norank2} 
Suppose that $C$ is a division algebra over $K$. Then there is a $14$-dimensional subspace of $\Alt(C_0)$ which  contains no elements of rank $2$ and is invariant under $\Aut(C)$.

\end{theorem}

It follows from Theorem 4 of \cite{GQ} that 14 is the largest possible dimension for a subspace of $\Alt(C_0)$ which  contains no elements of rank $2$. 

For a field such as $\mathbb{R}$, Corollary \ref{orbits} and Theorem \ref{norank2} admit interesting interpretations in terms of the fundamental representations of the compact group $G_2^c$. The action of $G_2^C$ on the subspace
$\M$ of $\Alt(C_0)$ may be identified with the fundamental representation of degree 7 of this group. Now $C_0\wedge C_0$
is the direct sum of two irreducible $G_2^c$-submodules, one of dimension 7 and the other of dimension 14. See, for example, the formula on p.353 of \cite{FH}. The 14--dimensional submodule affords the adjoint representation of the group, which is the other fundamental representation.

The subspace we have described in Theorem \ref{norank2} thus affords the 14-dimensional fundamental representation
of $G_2^c$. In particular, we can see that the two fundamental representations of this group may be realized
on subspaces of alternating bilinear forms on a seven-dimensional real vector space, one subspace containing
no elements of rank 4 or 6, the other no elements of rank 2. These are the largest dimensions for such subspaces, and their direct sum is the whole space of alternating bilinear forms.

Finally, when $C$ is a division algebra, we may imitate the proof of Corollary \ref{nonvanishing} to show that there is no two-dimensional subspace
of $C$ on which the elements of $\N$ all vanish. Then following the proof of Theorem \ref{kernel}, we obtain the following analogue of Theorem \ref{norank2}.

\begin{theorem} \label{norank2again} 
Suppose that $C$ is a division algebra over $K$. Then there is a $21$-dimensional subspace of $\Alt(C)$ which  contains no elements of rank $2$ and is invariant under $\Aut(C)$. Any subspace of $\Alt(C)$ of larger dimension  contains 
elements of rank $2$.

\end{theorem}


\begin{thebibliography}{9}


\bibitem{A} I. Agricola, \emph{Old and new on the exceptional group $G_2$}, Notices Amer. Math. Soc. $\mathbf{55}$ (2008), 922-929.
\smallskip
\bibitem{FH} W. Fulton and J. Harris, \emph{Representation theory}. Graduate Texts in Mathematics, 129, Springer-Verlag, New York, 1991.

\smallskip
\bibitem{GQ} R. Gow and R. Quinlan, \emph{On the vanishing of subspaces of alternating bilinear forms}, Linear Multilinear Algebra $\mathbf{54}$ (2006), 415-428.
\bibitem{LY} K. Y. Lam and P. Yiu, \emph{Linear spaces of real matrices of constant rank},  Linear Algebra Appl. $\mathbf{195}$ (1993), 69-79. 

\smallskip
\bibitem{SV} T. A. Springer and F. D. Veldkamp, \emph{Octonions, Jordan algebras and exceptional groups}, Springer-Verlag, Berlin, 2000.

\end{thebibliography}
\end{document}